\documentclass{amsart}
\usepackage{amssymb}  
\usepackage{amsmath}  
\usepackage{verbatim}
\usepackage{bm} 

\newcommand{\aplicacionsinelementos}[3]{#1 \stackrel{#2}{\longrightarrow }#3} 

        \newtheorem{theorem}{{\bf Theorem}}
        
        \newtheorem{lemma}[theorem]{{\bf Lemma}}
        
        \newtheorem{corollary}[theorem]{{\bf Corollary}}

        \newtheorem{definition}{{\bf Definition}} 
        \newtheorem*{theorem*}{\bf Theorem}
        \newtheorem*{lemma*}{\bf Lemma}
        \newtheorem*{proposition*}{\bf Proposition}
        \newtheorem*{corollary*}{\bf Corollary}
        \newtheorem*{remark*}{{\bf Remark}}
        \newtheorem*{example*}{{\bf Example}}
        \newtheorem*{definition*}{{\bf Definition}}
\theoremstyle{definition}
        \newtheorem{remark}[theorem]{{\bf Remark}}

\newcommand{\lh}{L'H\^{o}\-pi\-tal }
\newcommand{\lhs}{L'H\^{o}\-pi\-tal's }

\DeclareMathOperator{\ord}{ord}

\newcounter{problema}

\newcommand{\lhpunto}{L'H\^opital. }

\title{The Valuative Theory of Foliations}
\author{P. Fortuny Ayuso} 
\thanks{Partially supported by the European Commission: TMR Network ``Singularidades de Ecuaciones
  Diferenciales y Foliaciones'', ERBFMRXCT960040.\\
The author wishes to thank Prof. J. M. Aroca for his constant support during the preparation of his Ph. D. Thesis.}
\date{\today}
\begin{document}   
\vfuzz4pt 
\hfuzz3pt 
\begin{abstract}
  This paper gives a characterization of valuations which follow the singular infinitely near points of
  plane vector fields, using the notion of \lh valuation, which generalizes a well known classical
  condition. With that tool, we give a valuative description of vector fields with infinite
  solutions, singularities with rational quotient of eigenvalues in its linear part, and polynomial
  vector fields with transcendental solutions, among other results.  
\end{abstract}
\maketitle
\setcounter{section}{-1}
\section{Introduction}
The purpose of this paper is to characterize valuations which follow the infinitely near points given by a
singularity of a holomorphic vector field on a surface or, what is the same, a germ of a singular
first order differential equation in the
complex plane. Conceptually, a germ of
vector field  is nothing but a way of measuring infinitesimal approximations to a local family of
curves (the leaves of the associated foliation). Finding a solution through a point is the same as
obtaining  a subvariety which
has infinite contact with the field. Likewise, valuations are the generalization of contact along
a local object. Quoting Seidenberg \cite{Sei2}, ``(...) derivations are related to contact, and so are
valuations, so one may ask for a study connecting derivations and valuations''.\par
To that end, we use the concept of \lh valuation, which we define in Section $2$. It is a
generalization of what Rosenlicht calls \emph{differential valuations} in the context of Hardy fields
\cite{Ros}. A valuation $\nu$ is \lh for a vector field $\partial $ if it satisfies
the following generalization of classical \lhs condition:
\begin{equation*}
  \nu
  \left(
    \frac{a}{b}-\frac{\partial a}{\partial b}
  \right)>0
\end{equation*}
whenever $\nu(a)\geq\nu(b)>0$ and $\partial b\neq 0$. The key result of the paper is
Theorem $1$ of Section $3$, which states, essentially, that a \lh
valuation of a vector field follows (after a finite number of blowing-ups) infinitely near
singular points of the field. From this it follows
that valuations associated to curves which are not leaves of the foliation are not \lh for
the field. Section $3$ is devoted to the statement and proof of the main results. 
We single out two: the
equivalence between the existence of an infinite number of invariant analytic curves
through the origin (dicriticalness) and \lhs condition for divisorial valuations (Theorem
\ref{the:dicritic}), and the intrinsic
characterization of simple singularities with rational quotient of eigenvalues in terms of valuations of
rank $1$ and rational rank $2$ (Theorem \ref{the:generic-infinite}). In the last section,
we study the relation of rank $1$ valuations
with the existence of transcendental solutions, and a \emph{weak \lh condition}, which allows one to find
\emph{all} the solutions of a vector field, whether singular or not.\par
In the first section we recall, for the sake of both completeness and fixing of notation, the
classification of valuations centered in ${\mathbb C}\{x,y\}$ done by
Spivakovsky \cite{Sp} for $2-$dimensional regular local rings.\par
Let us remark also that the present work introduces a new algebraic perspective in the study of
singularities of plane and higher dimensional foliations (see \cite{tesis}), which up to
now have been studied  mainly
with complex analytic tools -either complex analysis or resolution of singularities.\par
This paper contains the research of the author's Ph. D. Thesis concerning complex plane
foliations. Results in higher dimension will appear in a forthcoming work.

\section{Valuations centered in ${\mathbb C}\{x,y\}$ and blowing-ups}
We are mainly interested in valuations centered in a point of the complex plane. Thus, we
shall work  over the ring ${\mathcal O}={\mathcal O}_{2}$ of germs of holomorphic
functions at a point $P\in{\mathbb C}^{2}$, near which we fix local analytic coordinates $(x,y)$, so
that ${\mathcal O}\simeq {\mathbb C}\{x,y\}$. More precisely, if ${\mathcal M}={\mathcal O}_{(0)}$
denotes the field of meromorphic functions at $P$, we are going to study valuations $\nu$ of ${\mathcal
  M}$, centered in the local ring $({\mathcal O},{\mathfrak m})$,
that is, such that $\nu({\mathcal O})\geq 0$ and $\nu({\mathfrak m})>0$. In this section we fix some notation
and terminology concerning these
objects. All the results are well-known since (at least) Zariski \cite{Zar4}, although our main
reference will be Spivakovsky's work \cite{Sp}, where one can find a general study of valuations centered 
in noetherian rings of
dimension $2$, without restrictions on the characteristic.
Let $\nu:{\mathcal M}^{\star}\rightarrow \Gamma$ be a valuation of ${\mathcal M}$ centered in ${\mathcal
  O}$ and consider the blowing-up $\pi :\mathcal{Y}\rightarrow ({\mathbb C}^{2},0)$ with
  center $(0,0)$. As $\pi$ is a birational morphism, the field of fractions of ${\mathcal
  Y}$ is ${\mathcal M}$, so that
it makes sense to speak of $\nu$ as {\em being centered at a point $Q\in{\mathcal Y}$}. The following result
links valuations in ${\mathcal M}$ centered in ${\mathcal O}$ with valuations centered at points of
$E\subset  {\mathcal Y}$. The proof is an easy exercise:
\begin{lemma*}
  \label{lem:blowing-up-valuations}
  One and only one of the following statements holds:
  \begin{enumerate}
  \item $\nu$ is the ${\mathfrak m}-$adic valuation at $({\mathbb C}^{2},0)$. That is, if $f\in
    {\mathcal O}={\mathbb C}\{x,y\}$, then $\nu(f)=\ord_{(x,y)}(f)$.
  \item There is one and only one closed point $Q\in E$ such that $\nu$ is centered in $Q$. Moreover, $\nu$ can be
    extended to a valuation  of $\tilde{{\mathcal M}}$, the field of fractions of
    $\tilde{{\mathcal O}}_{Q}$, which is centered in $\tilde{{\mathcal O}}_{Q}$. 
  \end{enumerate}
\end{lemma*}
(The extended valuation of $2$ will be denoted  by $\nu$, for the sake of simplicity).
This lemma is the basis of the classification of valuations made by Zariski \cite{Zar4} and
Spivakovsky \cite{Sp}, as it allows a grouping in terms of {\em co-final behaviour} under
blowing-ups. The first step is to distinguish between divisorial and non-divisorial valuations:
\begin{corollary*}
  \label{cor:chain-of-blowing-ups}
  Let $\nu$ be a valuation of ${\mathcal M}$ centered in ${\mathcal O}$. One and only one of the
  following statements holds:
  \begin{enumerate}
  \item There is a finite sequence of blowing-ups
    \begin{equation*}
      \aplicacionsinelementos{{\mathcal X}_{n+1}}{\pi _{n}}{{\mathcal X}_{n}}\rightarrow
      \cdots\aplicacionsinelementos{}{\pi _{1}}{{\mathcal X}_{1}}\aplicacionsinelementos{}{\pi
      _{0}}{{\mathcal 
        X}_{0}=({\mathbb C}^{2},0)}
    \end{equation*}
    (where $n\geq -1$)
    such that if $E_{i+1}$ is the exceptional line of $\pi_{i}$ and $Q_{i+1}\in E_{i+1}$ is the center
    of the blowing-up $\pi_{i+1}$, then $\nu$ is centered in $Q_{i+1}$ for any $i=-1,\dots,n-1$, and if
    $(u,v)$ are local coordinates at $Q_{n+1}\in{\mathcal X}_{n+1}$, then $\nu$ is the $(u,v)-$adic
    valuation. In this case, $\nu$ is called the {\em divisorial valuation with strict center
    $Q_{n+1}$}.
  \item There is an infinite sequence of blowing-ups
    \begin{equation*}
      \cdots\rightarrow \aplicacionsinelementos{{\mathcal X}_{n+1}}{\pi _{n}}{{\mathcal X}_{n}}\rightarrow
      \cdots\aplicacionsinelementos{}{\pi _{1}}{{\mathcal X}_{1}}\aplicacionsinelementos{}{\pi
      _{0}}{{\mathcal X}}
    \end{equation*}
    such that if $E_{i+1}$ is the exceptional line of $\pi_{i}$ and $Q_{i+1}\in E_{i+1}$ is the center
    of the blowing-up $\pi_{i+1}$, then $\nu$ is centered in $Q_{i+1}$ for any $i\geq -1$.
  \end{enumerate}
  In both cases, the chain of blowing-ups is said to be {\em associated to $\nu$}.
\end{corollary*}
In order to classify non-divisorial valuations, we introduce some terminology. Suppose $\nu$ is
non-divisorial and let ${\mathcal X}_{n+1}\stackrel{\pi}{\rightarrow }{\mathcal X}_{0}=({\mathbb
  C}^{2},0)$ be the composition of the first $n$ blowing-ups of the sequence
$(\pi_{i})_{i\geq 0}$ associated to $\nu$. The fiber $F_{n+1}$ of $0$ is the union of
$n+1$ irreducible components
\begin{equation*}
  F_{n+1}=E_{1}\cup \dots \cup E_{n+1},
\end{equation*}
where $E_{i+1}$ is the exceptional line of the blowing-up $\pi_{i}$.
\begin{definition*}
  The line $E_{i+1}$ is called the $(i+1)-$st irreducible component of $F_{n+1}$, for $i=0,\dots,n$.
\end{definition*}
\begin{definition*}
  We say that $Q_{n+1}$ is a {\em corner} if it belongs to two different irreducible components of
  $F_{n+1}$. Otherwise, we say that it is a {\em regular point} of $F_{n+1}$.
\end{definition*}
These notions allow one to  group non-divisorial valuations of ${\mathcal M}$ centered in ${\mathcal
  O}$ into four different families. Recall that any valuation $\nu$ of ${\mathcal M}$ can be naturally
  extended to one and only one valuation $\hat{\nu}$ of $\hat{{\mathcal M}}={\mathbb C}((x,y))$.
\begin{definition*}
   Fix a valuation $\nu$ centered in ${\mathcal O}$. We say that
  \begin{enumerate}
  \item $\nu$ is {\em associated to a curve}, or is {\em of contact with a curve} if there is an $n_{0}\in
    {\mathbb N}$ such that for any $m\geq n_{0}$, the point $Q_{m}$ is a regular point of $F_{m}$. In this
    case, there exists
    an irreducible principal ideal
    $(\hat{f})\in \hat{{\mathcal O}}={\mathbb C}[[x,y]]$ such that the valuation $\hat{\nu}$ has value
    group ${\mathbb Z}^{2}$, with the lexicographic order, and is given by
    the condition
    \begin{equation*}
      \nu(\hat{g})=(i,j)\Leftrightarrow \hat{g}\in(\hat{f})^{i}-(\hat{f})^{i+1} \text{ and }
      \#(\frac{\hat{g}}{\hat{f}},\hat{f})=j,
    \end{equation*}
    where $\#$ indicates the ordinary intersection multiplicity. There are two subcases: either
    $(\hat{f})\cap {\mathcal O}=(0)$, in which case $\nu$ is said to be 
    {\em associated to a formal    curve}, or $(\hat{f})\cap {\mathcal O}\neq(0)=(f)$, where
    $f$ is irreducible. In this case, $\nu$ is
    {\em associated to a convergent curve}. If no adjective is present, a valuation associated to a curve
    is assumed to be associated to a convergent one. Notice that in the non-convergent
    case, if $(t^{\alpha },\varphi (t))$ is a Puiseux parametrization of $(\hat{f}=0)$,
    then the value of any $g\in \mathcal{O}$, is exactly
    \begin{equation*}
      \nu(g)=\ord_{t}(g(t^{\alpha },\varphi (t))).
    \end{equation*}
  \item $\nu$ is {\em of contact with a divisor} if there is $n_{0}\in{\mathbb N}$ such that $m\geq n_{0}$
    implies $\pi_{m+1}$ is centered in the
    intersection of $E_{m}$ and $E_{n_{0}}$.  In
    fact, these valuations can be understood as associated to the curve $E_{n_{0}}$ at the
    point $E_{m}\cap E_{n_{0}}$ {\em in the birational model}
    ${\mathcal X}_{m}$ of ${\mathcal M}$. The points $Q_{m}=E_{n_{0}}\cap E_{m}$ will be
    called the \emph{fixed corners of $\nu$}. 
  \item $\nu$ {\em has an irrational Puiseux exponent} if for $m>>0$, the point $Q_{m}$ is a corner of
    $F_{m}$, but these corners do not share a common
    divisor (which is case 2): $\nu$ ``jumps'' from divisor to divisor, but is always centered in
    corners. For $m>>0$, there is a local system of coordinates $(u,v)$ at $Q_{m}$ and $\lambda \in
    {\mathbb R}_{>0}-{\mathbb Q}_{>0}$ such that if $f(u,v)=\sum a_{ij}u^{i}v^{j}\in \hat{{\mathcal
        O}}_{Q_{m}}\cap \hat{\mathcal{M}}$, then
    \begin{equation*}
      \nu(f)=\nu\left(
        \sum_{i,j>0}a_{ij}u^{i}v^{j}
      \right)=\min_{a_{ij}\neq 0}{\lambda i+j}.
    \end{equation*}
    This $\lambda $ depends on $m$, and can be constructed by a Bezout-type algorithm (see \cite{Sp}). 
  \item $\nu$ {\em has an infinite number of Puiseux exponents} if in the sequence of  centers of
    $\nu$ there  is
    an  infinite number of both corners
    and  regular points. In
    \cite{Sp}  one can see the proof of the following fact: 
    there is a minimal generating sequence $(Q_{i})_{i\in {\mathbb N}_{0}}$ (see \cite{Sp}) for $\nu$, such
    that each $Q_{i}$ defines an analytically irreducible curve $C_{i}$ with $i-1$ Puiseux exponents,
    and such that, for $0<i<j$, the curve $C_{i}$ has maximal contact with $C_{j}$,
    which explains the expression ``infinite Puiseux exponents''.
  \end{enumerate}
\end{definition*}
Valuations of contact with a divisor (type $2$) are to be considered a subclass of those of contact with
a curve, for up to birational morphisms, there is no difference between them. Thus, we shall only
mention them explicitly in Corollary \ref{cor:contact-with-divisor}.

\section{L'H\^opital Valuations}
Let ${\mathcal D}_{{\mathcal O}}={\mathcal D}er_{{\mathbb C}}^{{\mathfrak m}}({\mathcal O})$ denote the
${\mathcal 
  O}-$module of inner derivations of ${\mathbb C}\{x,y\}={\mathcal O}$ trivial over ${\mathbb C}$ and
{\em continuous} for the ${\mathfrak m}-$adic topology. In fact, ${\mathcal D}$ is the ${\mathcal
  O}-$module of germs of analytic vector fields at $({\mathbb C}^{2},0)$, which is free of rank
$2$. We shall denote by $\Omega_{{\mathcal O}}$  the module of germs of holomorphic $1-$forms over
$({\mathbb C}^{2},0)$, and by ${\mathcal D}_{{\mathcal M}}$ and $\Omega_{{\mathcal M}}$ the
${\mathcal M}-$vector spaces of meromorphic derivations and forms, respectively.
For reasons that will become
apparent later, we are interested in the projectivization of these vector spaces, namely
${\bm{\mathcal D}}={\mathcal D}-\{0\}/{\mathcal M}^{\star}$ and ${\bm {\Omega}}=\Omega-\{0\}/{\mathcal
  M}^{\star}$, which can be regarded (via the usual pairing) as the same set, called the set of 
{\em germs of foliations in $({\mathbb C}^{2},0)$}. A (germ of) {\em holomorphic foliation} will be an element
 $\bm{\partial}\in\bm{{\mathcal D}}$ (or equivalently, its dual $\bm{\omega }\in\bm{\Omega}$).\par
Let $\nu$ be a valuation of ${\mathcal M}$ and $\bm{\partial }\in\bm{\mathcal D}$ a foliation. Recall
that  one defines $\nu(0)=\infty$, where
$\infty$ is included in $\Gamma$, with the usual order and addition rules. The relation between foliations
and valuations is given by the following definition, which is a generalization of what Rosenlicht calls,
in a more restrictive context,
{\em differential valuations} \cite{Ros}:
\begin{definition}
  \label{def:lhopital}
  We say that $\nu$ is a {\em L'H\^opital valuation for} $\bm{\partial }$ if for any representative
  $\partial $ of $\bm{\partial }$, the following four equivalent conditions hold:
  \begin{enumerate}
  \item For any $a,b\in{\mathcal M}^{\star}$ with $\nu(a)\geq\nu(b)>0$ and $\partial b\neq 0$, 
    \begin{equation*}
      \nu\left(
        \frac{a}{b}-\frac{\partial a}{\partial b}
      \right)>0.
    \end{equation*}
  \item For any $a,b\in{\mathcal M}^{\star}$ with $0>\nu(a)\geq\nu(b)$ and $\partial b\neq 0$, 
    \begin{equation*}
      \nu\left(
        \frac{a}{b}-\frac{\partial a}{\partial b}
      \right)>0.
    \end{equation*}
  \item For any $a,b\in{\mathcal M}^{\star}$ with $\nu(a)\geq0$, $\nu(b)>0$ and $\partial b\neq 0$,
    \begin{equation*}
      \nu\left(
        \frac{\partial ab}{\partial b}
      \right)>0.
    \end{equation*}
  \item For any $a,b\in{\mathcal M}^{\star}$ with $\nu(a)\geq0$, $\nu(b)<0$ and $\partial b\neq 0$,
    \begin{equation*}
      \nu\left(
        \frac{\partial ab}{\partial b}
      \right)>0.
    \end{equation*}
  \end{enumerate}
\end{definition}
\begin{remark*}
  Notice that the definition is independent of the representative of $\bm{\partial }$, as the derivation
  appears both in the numerator and the denominator of the left-hand member of each equation.
\end{remark*}
The above notion is the natural translation of classical \lhs rule to germs of holomorphic vector fields
in the plane: the main difference lying, {\em grosso modo}, in the fact that in the one-variable case,
there are only one 
foliation and one valuation (compatible with the local structure of ${\mathbb C}\{t\}$), whereas in
higher dimensions one needs to fix the foliation and impose conditions on the valuation to relate both
concepts. 
\begin{example*}{See \cite{tesis} and \cite{Ros}}
  Consider the differential equation
  \begin{equation*}
    \frac{dy}{dx}=-\frac{y}{x^{2}},
  \end{equation*}
  having the fundamental solution $y=e^{1/x}$, which is holomorphic in a punctured neighbourhood of
  $0\in{\mathbb C}$. Let ${\mathcal K}={\mathbb C}(z,e^{1/z})$ be the field of rational
  functions in $z,e^{1/z}$. Fix two analytic paths $\gamma _{0},\gamma _{1}$, with $\gamma _{0}\subset
  \{\text{Re}(z)>0\}$, $\gamma _{1}\subset\{\text{Re}(z)<0\}$ and $\gamma _{0}(1)=\gamma _{1}(1)=0$.
  Define
  ${\mathcal O}^{i}=\{f\in{\mathcal K}:\lim_{t\rightarrow 1}\vert f(\gamma _{i}(t))\vert <\infty\}$, for
  $i=0,1$. Both ${\mathcal O}^{i}$ are valuation rings. Call $\nu_{i}$ to the
  valuation associated to ${\mathcal O}^{i}$, and let $\partial \in\text{Der}_{{\mathbb C}}{\mathcal K}$
  be the derivation of ${\mathcal K}$ induced by the differential equation: $\partial z=1,\partial
  e^{1/z}=-1/z^{2}e^{1/z}$. It is
  easy to check (see \cite{Ros}) that each $\nu_{i}$ is a \lh valuation for $\partial $. We remark that
  these rings depend only on ``the side of ${\mathbb C}$ where $\gamma $ is in'', that is, paths included
  in $\text{Re}(z)>0$ give the same ring ${\mathcal O}^{0}$ and paths in $\text{Re(z)}<0$ give ${\mathcal
  O}^{1}$. This leads us to think that \lh valuations are also a means to find Stokes lines for ordinary
  differential equations, although this remains, up to date, a conjecture. 
\end{example*}

\section{Valuations and vector fields in dimension $2$}\label{sec:valuations-fields-dim2}
We have presented the necessary background to study the links between germs of foliations  in $({\mathbb
  C}^{2},0)$ and valuations centered in ${\mathcal O}$. As we have remarked, we are going to use both
  the languages of $1-$forms and of vector fields (distributions), as they are the same thing
  in this dimension.\par
Take a holomorphic $1-$form $\omega $ and fix a regular system of parameters
  $(x,y)$ in $({\mathbb C}^{2},0)$. As ${\mathcal O}\simeq {\mathbb C}\{x,y\}$, $\omega $ can be
  written
\begin{equation*}
  \omega =adx+bdy,
\end{equation*}
with $a,b$ holomorphic functions in $({\mathbb C}^{2},0)$. Given such a form, we have its
line of dual vector
fields $[\partial]$ generated by
\begin{equation*}
  \left\{
    \begin{array}{l}
      \partial x=-b\\ \partial y =a
    \end{array}
  \right.
\end{equation*}
We might get any non-zero vector field in $[\partial ]$, by the remark following
Definition \ref{def:lhopital}. From now on, $\omega $ and the above $\partial$  are
fixed. We say that $(0,0)$ is a \emph{singular point} for a foliation $\bm{\omega}\in\bm{\Omega}$ if there is a
coordinate system $(x,y)$ such that for any representative $\omega=a(x,y)dx+b(x.y)dy$ of $\bm{\omega}$,
$a(0,0)=b(0,0)=0$.\par
By a {\em separatrix} of $\omega $ we shall mean a formal irreducible curve
tangent to $\omega $. That is, a separatrix is a principal irreducible ideal
$(\hat{f})\subset {\mathbb C}[[x,y]]$ such that $d\hat{f}\wedge \omega=\hat{f}\eta$, with $\eta$ a formal
$2-$form. From the  point of view of vector fields, a separatrix can be defined as a local non-zero ${\mathbb
C}-$morphism
\begin{equation*}
  \varphi :{\mathbb C}[[x,y]]\rightarrow {\mathbb C}[[t]],
\end{equation*}
such that
\begin{equation*}
  \left\{
    \begin{array}{l}
      {\partial \varphi (x)}/{\partial t}=-b(\varphi (x),\varphi (y))\\ 
      {\partial\varphi (y)}/{\partial t}=a(\varphi (x),\varphi (y))
    \end{array}
  \right.
\end{equation*}
\subsection{The key results}
Fix a valuation  $\nu$ of 
${\mathcal M}={\mathbb C}\{\{x,y\}\}$, $\nu$ centered in ${\mathcal O}$. Let ${\mathcal
Y}\stackrel{\pi}{\rightarrow }{\mathcal X}$ be the blowing-up of ${\mathcal X}=({\mathbb C}^{2},0)$ at
$0$ and call $E=\pi^{-1}(0)$, the exceptional divisor. Suppose there is a point $Q\in E$ such that $\nu$
is centered in $Q$ (that is, $\nu$ is not the ${\mathfrak m}-$adic valuation). Let ${\mathcal O}_{Q}$ be
the local ring at $Q$ and $\tilde{\omega }$ (respectively $\tilde{\partial}$) be the strict transform of
$\omega $ at $Q$ (resp. the strict transform of $\partial$). The following result is the cornerstone
of this work:
\begin{theorem}\label{the:cornerstone}
 If $\nu$ is a \lh
 valuation for $\partial$ and $E$ is invariant for $\tilde{\omega }$ (that is to say,
 $(0,0)$ is not dicritical for $\omega $), then $Q$ is a singular point for $\tilde{\partial}$.
\end{theorem}
\begin{proof}
  Suppose, in order to
  get a contradiction, that $Q$ is not singular for $\tilde{\partial}$. After a linear change of
  coordinates at $({\mathbb C}^{2},0)$, we may assume that $Q$ is the origin of the following chart of
  $\pi$:
  \begin{equation*}
    \pi: \left\{
      \begin{array}{l}
        x=x_{2}\\ 
        y=x_{2}y_{2}
      \end{array}
    \right.
  \end{equation*}
  Let $\omega$ be the dual form of $\partial $ and $\tilde{\omega }$ its strict transform. As $(0,0)$ is
  non-dicritical for $\omega$ and $Q$ is not singular for $\tilde{\omega }$, the starting
  form can be written as
  \begin{equation*}
    \omega =x_{2}^{\alpha }u(dx_{2}+x_{2}\psi dy_{2}),
  \end{equation*}
  where $\psi\in {\mathcal O}_{Q}$ (in fact, it is in ${\mathcal O}\subset {\mathcal O}_{Q}$) and $u$
  is an element of ${\mathcal O}$ not zero at $(x_{2}=0,y_{2}=0)$. The factor $x_{2}$ appears because
   $E$ is invariant ($\omega$ is non-dicritical at $(0,0)$). As $\nu(x_{2}),\nu(y_{2})>0$,
  it is clear that $\nu(y)>\nu(x)>0$. The vector field $\tilde{\partial}$ is then
  \begin{equation*}
    \left\{
      \begin{array}{l}
        \tilde{\partial}x_{2}=-x_{2}^{\alpha +1}u\psi\\ \tilde{\partial}y_{2}=x_{2}^{\alpha }u
      \end{array}
    \right.
  \end{equation*}
  whence
  \begin{equation*}
    \nu\left(\frac{y}{x}-\frac{\partial{y}}{\partial{x}}\right)=\nu\left( y_{2}-
      \frac{x_{2}\/\tilde{\partial}{y_{2}}+y_{2}\/\tilde{\partial}{x_{2}}}{\tilde{\partial}{x_{2}}}\right)=
  \end{equation*}
  \begin{equation*}
    =\nu\left( -\frac{x_{2}^{\alpha+1}u}{-x_{2}^{\alpha+1}u\psi} \right),
  \end{equation*}
  but $\nu(u)=0$, $\nu(1)=0$, $\nu(\psi)\geq 0$ and $\nu$ cannot be a \lh valuation for $\partial$.
\end{proof}
From Theorem \ref{the:cornerstone} and Seidenberg's reduction of singularities for plane
vector fields \cite{Sei}, we
infer immediately the following
\begin{corollary}\label{cor:infinite-pairs}
  Valuations with an infinite number of Puiseux exponents are never \lh valuations for any
  holomorphic vector field $\partial$.
\end{corollary}
\begin{proof}
  After a finite number of blowing-ups, we may assume that $\partial $ is non-dicritical.
  Suppose, in order to get a contradiction, that $\nu$ is a \lh valuation following both regular points and
  corners of the exceptional divisor in the long run. Seidenberg's result implies \cite{Sei}
  that from some $Q_{0}$ onwards, all these points are \emph{simple} for the strict
  transform of $\partial $. Assume, without loss of generality, that $Q_{0}$ is corner and
  both irreducible components $E$ and $E'$ of the divisor passing through it are separatrices.
  From this, we deduce that the only
  infinitely near singular
  points of $\partial $ at $Q_{0}$ after  the blowing-up $\pi$  with center $Q_{0}$ are the two
  corners of $\pi^{-1}(Q_{0})$ corresponding to $E$ and $E'$. Theorem
  \ref{the:cornerstone} implies that, from $Q_{0}$
  onwards, all the infinitely near points of $\nu$ are
  corners, whence $\nu$ cannot have an infinite number of Puiseux exponents.
\end{proof}
\begin{lemma}\label{lem:chain-rule}
  Let $\partial $ be a holomorphic vector field and $\hat{f}\in {\mathbb C}[[x,y]]$ a strictly formal
  (which means $(\hat{f})\cap {\mathcal O}=(0)$)
  separatrix  of $\partial $
  not tangent to $(x=0)$. Write $\partial =
  h\partial /\partial
  x+g\partial /\partial y$ with $(g,h)=1$ and take $u=a/b\in {\mathcal M}$. Suppose $(t^{\alpha },\varphi
  (t))$ is a  Puiseux
  parametrization of $\hat{f}=0$, with $\ord_{t}\varphi (t)>\alpha $. Then
  \begin{equation}
    \label{eq:chain-rule}
    \partial u(t^{\alpha },\varphi (t))= \frac{1}{\alpha
      t^{\alpha-1}}h(t)\frac{d}{dt}(u(t^{\alpha},\varphi(t))).
  \end{equation}
\end{lemma}
\begin{proof}
  We have, by the chain-rule
  \begin{equation*}
    \partial u=u_{x}\partial x +u_{y}\partial y=u_{x}h+u_{y}g=h(u_{x}+\frac{g}{h}u_{y}).
  \end{equation*}
  As $\hat{f}=0$ is strictly formal, $h(t^{\alpha },\varphi (t))\neq 0$ and, from the definition of
  separatrix,
  we get:
  \begin{equation*}
    \frac{g}{h}(t^{\alpha },\varphi(t))=\frac{\varphi'(t)}{\alpha t^{\alpha -1}},
  \end{equation*}
  so that, substituting this expression in the previous one:
  \begin{equation*}
    \partial(u)(t^{\alpha },\varphi(t))= \frac{1}{\alpha
    t^{\alpha-1}}h(t)\frac{d}{dt}(u(t^{\alpha},\varphi(t))),
  \end{equation*}
  as desired.
\end{proof}
\subsection{The dicritical case}
The link between valuations and dicritical centers is founded on the following
\begin{lemma}\label{lem:dicritic-previous}
  Let $\partial$ be a holomorphic vector field in $({\mathbb C}^{2},0)$ and let $\nu$ be the ${\mathfrak
    m}-$adic valuation centered at ${\mathcal O}$ (that is, $\nu(a)=n\in {\mathbb N}$ if and only if
    $a\in {\mathfrak m}^{n}-{\mathfrak m}^{n+1}$).
    Put $r=\min\{\nu(\partial x), \nu(\partial y)\}$. Then
  \begin{enumerate}
  \item $\nu(a\partial b - b\partial a)\geq \nu(a)+\nu(b)+r-1$, for any $a,b\in {\mathcal O}$.
  \item $\partial$ is dicritical in $(0,0)$ if and only if there are $a,b\in {\mathcal O}$ for which the
strict inequality holds. Moreover, in this case, $\nu(a)=\nu(b)$.
  \end{enumerate}
\end{lemma}
\begin{proof}
  Take $a,b\in {\mathcal O}\subset {\mathbb C}[[x,y]]$ with $\nu(a)=m$, $\nu(b)=n$. Call
  $\theta=\partial x$ and $ \eta=\partial y$ and write $a,b,\theta$
  and $\eta$ as power series expansions:
  \begin{equation*}
    \begin{array}{ll}
      \left\{
        \begin{array}{l}
          a=a_{m}+a_{m+1}+\dots\\ b=b_{n}+b_{n+1}\dots
        \end{array}
      \right. & \left\{
        \begin{array}{l}
          \theta=\theta_{l}+\theta_{l+1}+\dots\\ \eta=\eta_{k}+\eta_{k+1}\dots
        \end{array}
      \right.
    \end{array}
  \end{equation*}
Applying $\partial$ to $a$ and $b$, we get
\begin{equation*}
  \left\{
\begin{array}{l}
\partial a=a_{mx}\theta_{l}+a_{my}\eta_{k}+\dots\\ \partial b=b_{nx}\theta_{l}+b_{my}\eta_{k}+\dots
\end{array}
\right.
\end{equation*}
where subindices $x$ and $y$ mean ``ordinary partial derivation''. From this last expression, it follows that
\begin{equation*}
  \nu(a\partial b-b\partial a)\geq m+n-1+r,
\end{equation*}
which is the general inequality.\par
To prove the second assertion, notice that when applying $\partial$
to a monomial $a_{ij}x^{i}y^{j}$, one gets
\begin{equation*}
  \partial (a_{ij}x^{i}y^{j})=i a_{ij}x^{i-1}y^{j}\theta + ja_{ij}x^{i}y^{j-1}\eta
\end{equation*}
and, looking at the terms of least degree in $\theta$ and $\eta$, we have
\begin{equation*}
  \partial (a_{ij}x^{i}y^{j})=i a_{ij}x^{i-1}y^{j}\theta_{r}+ja_{ij}x^{i}y^{j-1}\eta_{r}+\text{ h.o.t. }
\end{equation*}
hence
\begin{align*}
\partial (a_{ij}x^{i}y^{j})b_{pq}x^{p}y^{q}& -\partial (b_{pq}x^{p}y^{q})a_{ij}x^{i}y^{j}=\\
&=(ia_{ij}x^{i-1}y^{j}\theta_{r}+ja_{ij}x^{i}y^{j-1}\eta_{r})b_{pq}x^{p}y^{q}-\\ & -(p
b_{pq}x^{p-1}y^{q}\theta_{r}+qb_{pq}x^{p}y^{q-1}\eta_{r})a_{ij}x^{i}y^{j}+\text{ h.o.t }=\\
=((i-p)a_{ij}b_{pq}\theta_{r})&x^{i+p-1}y^{j+q}+ ((j-q)a_{ij}b_{pq}\eta_{r})x^{i+p}y^{j+q-1}+\text{
h.o.t. }
\end{align*}
which is zero if and only if $i+j=p+q$ (that is, $\nu(a)=\nu(b)$) and $x\theta_{r}+y\eta_{r}=0$ (i. e.
$\partial$ is dicritical at $(0,0)$), and the conclusion follows.
\end{proof}
Lemma \ref{lem:dicritic-previous} allows us to classify dicritical vector fields in terms of \lh valuations:
\begin{theorem}\label{the:dicritic}
  Let $\nu$ be a divisorial valuation of ${\mathcal M}$, centered at ${\mathcal O}$. Let $Q$ be the
  {\em strict center} of $\nu$ and $\partial$ a foliation by curves on $({\mathbb C}^{2},0)$. Then
\begin{center}
  $\nu$ is a \lh valuation for $\partial$ $\Leftrightarrow$ $\partial$ is dicritical at $Q$.
\end{center}
\end{theorem}
\begin{proof}
  Fix a local system of parameters $(u,v)$ at $Q$, $u,v\in {\mathcal O}_{Q}$. There are $\sigma, \theta,
  \eta\in {\mathcal O}_{Q}$ with coefficients in ${\mathcal
  O}_{Q}$ such that $\tilde{\partial}$ can be written as
  \begin{equation*}
    \tilde{\partial }=\frac{\theta(u,v)}{\sigma(u,v)}\frac{\partial}{\partial
    u}+\frac{\eta(u,v)}{\sigma(u,v)}\frac{\partial}{\partial v}.
  \end{equation*}
  Moreover, in those coordinates, $\nu=\ord_{(u,v)}$ is ``the order in $(u,v)$''. We can forget the
  denominators and suppose $\theta$ and $\eta$ are relatively prime in ${\mathcal O}_{Q}$. Let $m$ be the
  smallest among $\ord_{(u,v)}(\theta)$, $\ord_{(u,v)}(\eta)$. For any two meromorphic functions
  $a=\overline{a}(u,v)/\Delta(u,v)$, 
  $b=\overline{b}(u,v)/\Delta(u,v)$ with
  $\overline{a},\overline{b},\Delta\in {\mathcal O}_{Q}$, one has
  \begin{equation*}
    \frac{a}{b}-\frac{\partial {a}}{\partial {b}}=\frac{\Delta(\overline{a}\partial
      {\overline{b}}-\overline{b}\partial {\overline{a}})} {\overline{b}(\Delta\partial
      {\overline{b}}-\overline{b}\partial {\Delta})}.
  \end{equation*}
  Let $d=\ord(\Delta)$, $s=\ord(\overline{b})$, $t=\ord(\overline{a})$ and suppose, in
  order to study \lhs condition, that $\nu(a)\geq \nu(b)>0$, which means that $0\leq d<s\leq
  t$. Lemma \ref{lem:dicritic-previous} asserts that
  \begin{equation*}
    \nu(\Delta(\overline{a}\partial {\overline{b}}-\partial {\overline{a}}\overline{b}))\geq
    {\overline{b}(\Delta\partial {\overline{b}}-\overline{b}\partial {\Delta})}=2s+d+m-1.
  \end{equation*}
  The last equality hols because $\ord(\overline{b})>\ord(\Delta)$. \lhs condition means that the order of
  the numerator is greater that the order of the denominator. This follows easily if
  $t>s$. When $t=s$, the condition follows if and only if (cf. Lemma \ref{lem:dicritic-previous})
  $\partial$ is dicritical at $Q$, which finishes the proof.
\end{proof}
\subsection{Separatrices}
The first result relating solutions of differential equations with \lh valuations deals with strictly formal
separatrices:
\begin{theorem}\label{the:formal-separatrices}
  Let $\partial$ be a foliation by curves in $({\mathbb C}^{2},0)$ and $\nu_{\hat{f}}$ the valuation
  associated to a formal curve $(\hat{f}=0)$. Then
  \begin{center}
    $\hat{f}=0$ is a separatrix of $\partial$ $\Leftrightarrow$ $\nu$ is a \lh valuation for $\partial$.
  \end{center}
\end{theorem}
\begin{proof}
  $\Leftarrow)$ This is an easy consequence of Theorem \ref{the:cornerstone} and the fact that a (formal)
  curve following singular points of the strict transforms of $\partial$ by  blowing-ups is a
  separatrix.\par
  $\Rightarrow)$ $\nu$ has rank $1$ and, in fact, is the intersection multiplicity with
  $\hat{f}$, which means that if $(x=t^{\alpha },y=\varphi (t))$ is a Puiseux parametrization of
  $\hat{f}=0$ with $\ord_{t}(\varphi(t))>\alpha $ then, for $a=a(x,y)\in {\mathcal O}$, we have
  $\nu(a)=\ord_{t}(a(t^{\alpha },\varphi (t)))$. Let us prove that, in fact, it is a \lh valuation with
  respect to $\partial$. Write $\partial=h\frac{\partial}{\partial x}+g\frac{\partial}{\partial y}$ with
  $h,g\in {\mathcal O}$ and relatively prime in ${\mathcal O}$. By Lemma \ref{lem:chain-rule}, if $p,q\in
  {\mathcal M}$, then
  \begin{equation*}
    \nu \left( \frac{p}{q}-\frac{\partial{p}}{\partial{q}}\right)=ord_t \left(\frac{p}{q}(t)-
      \frac{\partial{p}}{\partial{q}}(t)\right) = ord_t \left(\frac{p(t)}{q(t)}-
      \frac{{p'}(t)}{{q'}(t)}\right) >0,
  \end{equation*}
  the last equality holding as a consequence of \lhs rule for one complex variable.
\end{proof}
It appears that in order to study valuations associated to convergent separatrices, one needs to distinguish
them according to the quotient of eigenvalues of the linear part of the vector field at a point. To be
more precise, let $\hat{f}$ be
a separatrix,  convergent or not, of a holomorphic vector field $\partial $ and let $\rho $ be a
resolution  of singularities
of $\hat{f}$. Put $E=\rho^{-1}(0,0)$. We may assume that the last exceptional line appearing in $E$ is
invariant for $\rho ^{\star}\partial $, that the strict transform of $\hat{f}$ meets $E$ in a regular
point of $E$ and, finally, that $Q$ is a simple point for $\rho ^{\star}\partial $ (cf. \cite{Sei}). Let
$(u,v)$ be a local regular system of parameters of ${\mathcal O}_{Q}$ such that $u=0$ is the equation of
$E$ and $v=0$ is tangent to the strict transform of $\hat{f}$. In these coordinates, the equation of
$\rho^{\star}\partial $ has non-zero linear part of the form $\mu u\partial /\partial u-\lambda v\partial
/\partial v$.
\begin{definition}
  We say that $\hat{f}$ has
  \begin{itemize}
  \item {\em Generic character} if $\lambda \mu \neq 0$ and $\lambda /\mu \not\in {\mathbb Q}$.
  \item {\em Infinite character} if $\mu =0$.
  \item {\em Zero character} if $\lambda =0$.
  \item {\em Rational character} if $\lambda \mu \neq 0$ and $\lambda /\mu \in {\mathbb Q}$.
  \end{itemize}
\end{definition}
\begin{remark}
  The transformation rule for of the coefficients of the linear part of $\partial $ by blowing-up shows that the
  character is well defined: it does not depend either on $\rho $ or on the coordinates $(u,v)$.
\end{remark}
\begin{theorem}
  \label{the:generic-infinite}
  If $\nu$ is of contact with the convergent curve defined by the ideal $(f)\subset{\mathcal O}$, then 
  \begin{center}
    $\nu$ is a \lh valuation $\Leftrightarrow$ $f$ is a separatrix with $\left\{
      \begin{array}{l}
        \text{ generic }\\ \text{ or }\\ \text{ infinite }
      \end{array}
    \right\}$ character.
  \end{center}
\end{theorem}
\begin{proof}
  For $a\in \mathcal{O}$, $\nu$ is given by
  \begin{equation*}
    \nu(a)=(i,j)\Leftrightarrow a\in(f)^{i}-(f)^{i+1},\,\#(\frac{a}{f^{i}},f)=j,
  \end{equation*}
  where $\#(a,b)$ is, as usual, the intersection multiplicity of the curves $(a=0)$ and $(b=0)$.\par
  $\Rightarrow$) From Theorem \ref{the:cornerstone}, we infer that $(f=0)$ must be a separatrix of
  $\partial $. As for the character, we divide the proof in two parts:\par
  {\bf a)} Assume, first, that
  $(f=0)$ has {\em rational character}. After reducing the singularities of $\partial $ (and hence of
  $f$), we 
  may assume that the eigenvalues of its linear part are $1$ and $m/n$, with $m,n\in {\mathbb N}$. We
  may also assume that the irreducible component $F$ of the exceptional divisor $E$ meeting $(f=0)$ is
  both invariant for $\partial $ and transversal to $(f=0)$ at a regular point of $E$. Take a local
  system of parameters at $Q$, say $(u,v)$, such that the strict transform of $f$ follows the same
  infinitely near points as $(u=0)$ and such that $F$ is given by $(v=0)$. Things being so, the vector field
  can be written, at $Q$, as follows:
  \begin{equation*}
    \left\{
      \begin{array}{l}
        \partial u=u(1+f_{1})\\ 
        \partial v=-\frac{m}{n}v(1+f_{2})
      \end{array}
    \right.\, \text{ with }\nu(f_{1}),\nu(f_{2})>0,
  \end{equation*}
  as $f=0$ and $v=0$ are invariant and $\partial $ is simple at $Q$. Take
  \begin{equation*}
    \left\{
      \begin{array}{l}
        a=u^{m}v^{n+1}\\ b=u^{m}v^{n}
      \end{array}
    \right.
  \end{equation*}
  By construction, $\nu(a)\geq\nu(b)>0$. Let us see if they satisfy \lhs criterion:
  \begin{equation*}
    \frac{a}{b}-\frac{\partial a }{\partial b }= \frac{v^{2n+1}(m/n+g_1 )}{v^{2n}g_2 }
    =\frac{\alpha}{\beta}
  \end{equation*}
  with $\nu(g_{1}),\nu(g_{2})>0$. As $(n+1)m/n\neq 0$, we have $\nu(\alpha )=(0,2n+1)$ and $\nu(\beta
  )>(0,2n)$, so that the pair $a,b$ does not satisfy the criterion and $\nu$ is not \lh for
  $\partial$.\par
{\bf b)} Suppose that $f=0$ has {\em character zero}. The same argument as in {\bf a)} shows
  that we can take a local system of coordinates $(u,v)$ at $Q$ with
  \begin{equation*}
    \left\{
\begin{array}{l}
\partial u=u^{q+1}+uf_{1}\\ \partial v=-v(1+f_{2})
\end{array}
\right.\, \text{ with } \nu(f_{1}),\nu(f_{2})>0,q\geq 1.
  \end{equation*}
  In this situation, take
  \begin{equation*}
    \left\{
\begin{array}{l}
a=uv\\ b=u(u+v)/v
\end{array}
\right.
  \end{equation*}
  which satisfy $\nu(a)\geq \nu(b)>0$. An easy computation shows that
  \begin{equation*}
    \frac{a}{b}-\frac{\partial a}{\partial b}=\frac{v^2}{u+v}- \frac{-v^2 -v^2 f_2+v^2 f_1+u^q v^2}{vf_1
    +ug},
  \end{equation*}
  where $g\in{\mathcal O}_{Q}$. As $\nu(u+v)=\nu(v)=(0,1)$ and $\nu(f_{1})>0$, the valuation does not
  satisfy \lhs criterion, q.e.d.\par
$\Leftarrow$) Fix, as before, $Q$ in the exceptional divisor of a
  reduction of singularities of $\partial $ and $(f=0)$, which is not a corner and take 
  coordinates $(u,v)$ at $Q$ such that $\partial $ is given by 
  \begin{equation*}
    \left\{
      \begin{array}{l}
        \partial u=u(1+f_{1})\\ \partial v=-\lambda v-vf_{2}
      \end{array}
    \right.  \,\text{ with } \nu(f_{1}),\nu(f_{2})>0.
  \end{equation*}
  Let $a=u^{i}\overline{a}$, $b=u^{j}\overline{b}$ with $\nu(a)\geq \nu(b)>0$ and
  $\nu(\overline{a})=(0,k)$, $\nu(b)=(0,l)$. We have
  \begin{equation}\label{eq:alpha-beta}
    \frac{a}{b}-\frac{\partial a }{\partial b }= \frac{u^{i+j}}{u^{2j}}\left( \frac
      {j\overline{a}\overline{b}+j\overline{a}\overline{b}f_1 +\overline{a}\partial{\overline{b}}-
      i\overline{a}\overline{b} -i\overline{a}\overline{b}f_1-\partial{\overline{a}}\overline{b}}
      {\overline{b}((j\overline{b}+\partial{\overline{b}})+j\overline{b}f_1)}
      \right)=u^{i-j}\frac{\alpha}{\beta}.
  \end{equation}
  Let us make a brief digression in order to calculate $\nu(\beta )$.\par
  \begin{minipage}[c]{348pt}
    \hspace{-5pt}\line(0,1){100}\line(0,-1){2350}
    \vspace{-235pt}
    \hspace{5pt}
    Write $\overline{b}=p/q$,  with $p,q\in {\mathcal O}_{Q}$ as follows:
    \begin{equation*}
      \begin{array}{l}
        p=v^{n}+uh_{1}+v^{n+1}h_{2}\\ 
        q=v^{m}+ug_{1}+v^{m+1}g_{2}
      \end{array}
    \end{equation*}
    where the functions $h_{i}$ and $g_{j}$ are formal power series in
    ${\mathbb C}[[u,v]]$. Their derivatives are
    \begin{equation}\label{eq:partials-pq}
      \begin{array}{l}
        \partial{p}=-n\lambda v^{n}+u\delta_1 -\lambda v^{n+1}\delta_2 +v^{n} f_2 \delta_3\\
        \partial{q}=-m\lambda v^{m}+u\epsilon_1 -\lambda v^{m+1}\epsilon_2 +v^{m}f_2 \epsilon_3
      \end{array}
    \end{equation}
    where, again, $\delta_{i}, \varepsilon _{j}$ are in ${\mathcal O}_{Q}$. From \eqref{eq:partials-pq},
    we have that $\nu(q\partial p-p\partial q)\geq 2\nu(q)$ and, hence, $\nu(\partial \overline{b})\geq
    \nu(\overline{b})$. In fact, equality holds if and only if $n\neq m$ and $\lambda \neq 0$. If
    $\nu(\partial \overline{b})>\nu(\overline{b})$, then $\nu(\beta )=2\nu(\overline{b})$. If not, then
    \begin{equation*}
      j\overline{b}+\partial{\overline{b}}=j\frac{v^{n}+\cdots}{v^{m}+\cdots}+\frac{\lambda(n-
        m)v^{n+m}+\cdots}{v^{2m}+\cdots}=\star,
    \end{equation*}
    (dots indicating terms of greater value). But
    \begin{equation*}
      \star=\frac{jv^{m+n}+\lambda(m-n)v^{n+m}+\cdots}{v^{2m}+\cdots}
    \end{equation*}
    and as, by hypothesis, $\lambda $ is neither zero nor rational, then $j+\lambda (m-n)\neq 0$ when $j\neq 0$,
    and in this case we get $\nu(j\overline{b}+\partial \overline{b})=\nu(\overline{b})$. If $j=0$, then
    $n>m$ and hence $\nu(\partial \overline{b})=\nu(\overline{b})$.\par
    Thus, in any case, $\nu(\beta )=2\nu(\overline{b})$.
\end{minipage}\par
  \noindent Continuing with the proof,  if $i>j$, then \lhs
  condition is satisfied because of the factor $u^{i-j}$. If $i=j$, then $\nu(\overline{a})\geq
  \nu(\overline{b})>0$, and from the expression of $\alpha $ in \eqref{eq:alpha-beta}, we infer that
  \begin{equation*}
    \nu (\alpha)\geq \nu
    \left(\overline{b}^{2}\partial\left(\frac{\overline{a}}{\overline{b}}\right)\right),
  \end{equation*}
  from what follows, by a reasoning similar to the previous one, that $\nu(\alpha )>\nu(\beta )$, and we
  are done.
\end{proof}
A straightforward consequence is the following
\begin{corollary}
  \label{cor:contact-with-divisor}
  Let $\nu$ be a valuation of contact with a divisor $E_{k}$, and let $\partial $ be a vector field. Let
  $Q$ be
  a fixed corner of $\nu$,  simple for $\partial $ and such that both components of the
  exceptional divisor passing through $Q$ are invariant. Then $\nu$ is a \lh valuation
  for $\partial $ if and only if one of the following conditions holds:
  \begin{itemize}
  \item[a)] The two irreducible components of the exceptional divisor passing through $Q$ have generic
    character (understood as separatrices in ${\mathcal O}_{Q}$).
  \item[b)] $E_{k}$ has infinite character (and hence, the other component has character zero).
  \end{itemize}
\end{corollary}

\subsection{Valuations with an irrational Puiseux exponent}
As in the previous case, we have to distinguish between rational and non-rational quotients of
eigenvalues, although in a slightly different way.
\begin{lemma}
  \label{lem:irrational-simple-nonrational}
  Let $Q$ be a simple singularity of a plane vector field $\partial =b_{1}\partial /\partial
  x-a_{1}\partial /\partial y$, where $a_{1},b_{1}$ are of order $1$ in $(x,y)$. Assume that the
  eigenvalues of the linear part of $\partial $ are both non-zero and have non-rational quotient. Then
  any valuation $\nu$ with an irrational Puiseux exponent following singular infinitely near points
  of $\partial $ is a \lh valuation for $\partial $.
\end{lemma}
\begin{proof}
  Up to a linear change of coordinates, we can write
  \begin{equation*}
    \left\{
      \begin{array}{l}
        \partial u=\lambda u+\dots\\
        \partial v=v+\dots
      \end{array}
    \right.
  \end{equation*}
  where $(u,v)$ is the new coordinate system. Also we may assume that the associated valuation is
  $\nu(u)=1$, $\nu(v)=\kappa \in {\mathbb I}_{+}$, so that if $\alpha =\sum a_{ij}u^{i}v^{j}\in {\mathbb
  C}[[u,v]]$, then $\nu(\alpha )=\min_{a_{ij}\neq 0}(i+\kappa j)=\ord_{t}(\sum a_{ij}t^{i+\kappa
  j})$. Given $\alpha ,\beta \in {\mathcal O}_{Q}$, we have
  \begin{equation}\label{eq:irrational1}
    \partial \left(\frac{\alpha}{\beta}\right)=\frac{(\lambda u\alpha_{u}+v\alpha_{v})\beta-(\lambda
      u\beta_{u}+v\beta_{v})\alpha}{\beta^{2}}+\ldots
  \end{equation}
  The set of monomials in $(\lambda u\alpha _{u}+v\alpha _{v})$ whose coefficient is non-zero is the same
  as that of $\alpha $, as $\lambda $ is  irrational. The same happens with $(\lambda u \beta_{u}+v\beta
  _{v}) $ and $\beta $. Let $(i_{0},j_{0})$ and $(m_{0},n_{0})$ be the leading coefficients of $\alpha $
  and $\beta$, respectively. The numerator in \eqref{eq:irrational1} is
  \begin{equation*}
    (\lambda i_0 +j_0 )\alpha_{i_0 j_0}\beta_{m_0 n_0}-(\lambda m_0 +n_0 )\alpha_{i_0 j_0}\beta_{m_0
    n_0}=
  \end{equation*}
  \begin{equation*}
    =\alpha_{i_0 j_0}\beta_{m_0 n_0}(\lambda (i_0 - m_0 )+(j_0 -n_0 )),
  \end{equation*}
  which is non-zero. If $(i_{0},j_{0})\neq (m_{0},n_{0})$, then the value of the numerator is the sum of
  the values of $\alpha $ and $\beta $. If, on the contrary, $(i_{0},j_{0})= (m_{0},n_{0})$, then the
  value of the numerator is strictly greater than the sum of the values of $\alpha $ and $\beta $.\par
  Take now $\alpha /\beta ,\gamma /\beta \in {\mathcal M}$. We need to study
  \begin{equation}\label{eq:numer-denom-irrational}
    \nu \left( \frac{\alpha /\beta }{\gamma /\beta }-\frac{\partial (\alpha /\beta )}{\partial (\gamma
    /\beta )} \right)= \nu \left( \frac{\beta (\alpha \partial \gamma -\gamma \partial \alpha )}{\gamma
    (\beta \partial \gamma -\gamma \partial \beta )} \right).
  \end{equation}
  Suppose $\nu(\alpha )\geq \nu(\gamma )>\nu (\beta )\geq 0$. From the second inequality and the previous
  considerations, it follows that $\nu(\gamma (\beta \partial \gamma-\gamma \partial \beta ))=\nu(\gamma
  )+\nu(\gamma )+\nu(\beta )$. There are two possibilities: either $\nu(\alpha )>\nu(\beta )$ or $\nu(\alpha
  )=\nu(\beta )$. In the first case, the value of the numerator in \eqref{eq:numer-denom-irrational} is,
  at least, $\nu(\alpha )+\nu(\beta )+\nu (\gamma )$, which is greater than the value of the
  denominator. To see that the same holds for the second case one only needs to use the above argument on
  the initial components of $\alpha ,\beta ,\gamma $ and their derivatives.
\end{proof}
\begin{remark*}
  Obviously, a valuation with an irrational Puiseux exponent which does not follow singular infinitely
  near points of a vector field is not a \lh valuation, in view of Seidenberg's reduction of
  singularities and Theorem \ref{the:cornerstone}.
\end{remark*}
\begin{lemma}
  \label{lem:irracional-simple-rational}
  Under the same conditions as in Lemma \ref{lem:irrational-simple-nonrational}, but with rational quotient of
  eigenvalues (including the case in which one of them is $0$), no valuation with an irrational
  Puiseux exponent is \lh for $\partial $.
\end{lemma}
\begin{proof}
  By the previous remark, we only need to consider valuations following singular infinitely near points
  of $\partial $. We divide the proof in two cases: both eigenvalues are non-zero (resonant case) and one
  of them is zero (saddle-node).\par 
  {\bf a) Resonant singularities.} We can write $\partial $ as
  \begin{equation*}
    \left\{
      \begin{array}{l}
        \partial u=-mu+\dots\\ \partial v=nv+\dots
      \end{array}
    \right.
  \end{equation*}
  After a polynomial change of coordinates (see \cite{CS}), we may assume that
  \begin{equation*}
    \partial =-(mu+s_{1}(u,v))\frac{\partial }{\partial u}+(nv+s_{2}(u,v))\frac{\partial}{\partial v},
  \end{equation*}
  with $\nu(s_{1}),\nu(s_{2})\geq \nu(u)+\nu(v)$ and that the valuation is $\nu(u)=1,\nu(v)=\eta\in
  {\mathbb I}_{>0}$ with $\eta>1$\footnote{Here we are using the fact that $m,n\neq 0$.}. Take
  $(i_{0},j_{0})\neq (0,0)$ such that $(-mi_{0}+nj_{0})=0$. For $cu^{i_{0}}v^{j_{0}}$, we have
  \begin{align*}
    \nu(\partial (cu^{i_{0}}&v^{j_{0}}))=\nu(i_0 cu^{i_0 -1}v^{j_0}s_1 (u,v)+ j_0 cu^{i_0}v^{j_0
      -1}s_2)\geq\\ 
    &\geq\min\{i_{0}+j_{0}\eta+\eta,i_{0}+j_{0}\eta+1\}=i_{0}+j_{0}+1,
  \end{align*}
  by definition. Let $r=[i_{0}+j_{0}\eta+1]$ (integral part). As $j_{0}\neq 0$ and $\eta\in {\mathbb
  I}$, we have $r<i_{0}+j_{0}\eta+1$. To verify that \lhs condition is not satisfied, take $\alpha
  =u^{r},\gamma =u^{i_{0}}v^{j_{0}}$, for which $\nu(\alpha )>\nu(\gamma )>0$, but $\nu(\partial \gamma
  )>\nu(\partial \alpha )=\nu(\alpha )$, so that
  \begin{equation*}
    \nu\left( \frac{\alpha }{\beta }-\frac{\partial \alpha }{\partial \beta } \right)=\nu(\alpha
)+\nu(\gamma )-\nu(\gamma )-\nu(\partial \gamma )<0,
  \end{equation*}
  which finishes the proof in this case.\par
  {\bf b) Saddle-node.} The vector field has a representative
  of the form
  \begin{equation*}
    \partial =a(u,v)\frac{\partial}{\partial u}+(v+b(u,v))\frac{\partial}{\partial v},
  \end{equation*}
  with $\ord_{(u,v)}(a),\ord_{(u,v)}(b)\geq 2$. The valuation is defined by $\nu(u)=1,\nu(v)=\eta\in
  {\mathbb I}_{>0}$, but we may not assume $\eta>1$. One can prove \cite{CS} that, up to a
  holomorphic change of coordinates, $\partial $ can be written
  \begin{equation*}
    \partial = (u(1+\lambda
    v^{p})+vR(u,v))\frac{\partial}{\partial u}+
    v^{p+1}\frac{\partial}{\partial v}
  \end{equation*}
  for suitable $\lambda \in {\mathbb C}$ and $p\in {\mathbb Z}_{\geq 1}$.
  Truncating, we deduce that there is $p\geq 1$ and a polynomial change of
  coordinates such that we may assume
  $\nu(a(u,v))\geq2$, (taking $a(u,v)=u^{p+1}+\dots$). Thus,
  for any $m\in{\mathbb N}$, we have
  \begin{equation*}
    \nu(\partial (u^{m}))\geq m+1.
  \end{equation*}
  Here we have to distinguish between $\eta<1$ and $\eta>1$. In the first case, let $\alpha =uv$ and
  $\gamma =u$. A direct computation shows that
  \begin{equation*}
    \nu\left( \frac{\alpha }{\gamma }-\frac{\partial \alpha }{\partial \gamma } \right)=\nu(\alpha
)+\nu(\gamma )-\nu(\gamma )-\nu(\partial \gamma )\leq 0.
  \end{equation*}
  In the second case, taking $r=[\eta]$ and $\alpha =v,\gamma =u^{r}$, the same argument finishes the
  proof.
\end{proof}
From these two lemmas, we infer the result concerning all the valuations with an irrational Puiseux
exponent:
\begin{theorem}
  \label{the:irrational-general}
  Let $\partial $ be a holomorphic vector field and $\nu$ a valuation centered in ${\mathcal O}$,
  with an irrational Puiseux exponent. Let $(P_{i})_{i\geq 0}$ be the sequence of infinitely near points of
  $\nu$. The following conditions are equivalent
  \begin{itemize}
  \item[a)] $\nu$ is a \lh valuation for $\partial $.
  \item[b)] $P_{i}$ is a singular point for the strict transform $\partial _{i}$ of $\partial $, for all
    $i$, and if $P_{k}$ is simple for $\partial _{k}$, then  both eigenvalues of its linear part are
    non-zero and their quotient is not in ${\mathbb Q}_{\leq 0}$.
  \end{itemize}
\end{theorem}
A consequence of the results of this section is the following characterization of valuations with an infinite
number of Puiseux exponents as those which are \emph{never} \lh valuations for any vector
field:
\begin{corollary}
  \label{cor:infinite-never}
  A valuation $\nu$ of $\mathcal{M}$ has an infinite number of Puiseux exponents if and only if there is
  no vector field for which it is \lhpunto
\end{corollary}
Sufficiency is Corollary \ref{cor:infinite-pairs}. The converse is derived from the
following lemma, whose prove is done by induction on the number of blowing-ups:
\begin{lemma}
  \label{lem:constructibility}
  Let $P$ denote an infinitely near point of $(0,0)\in {\mathbb C}^{2}$ appearing after a
  finite sequence of blowing-ups $\pi$. Then
  \begin{enumerate}
  \item If $P$ is a
    corner, then there exists a
    singular holomorphic foliation  $\omega $ in $({\mathbb C}^{2},0)$ such that the quotient of
    eigenvalues of the strict transform $\tilde{\omega }$ at $P$ is not a rational
    number. 
  \item If $P$ is not a corner  then, given any (formal or convergent) curve $\gamma $ passing
    through $P$ and transversal to the exceptional divisor, there exists a singular
    holomorphic foliation $\omega $ such that: $\gamma
    $ is a separatrix of $\omega $, $P$ is a simple singularity of the strict transform
    $\tilde{\omega }$ of $\omega $, and the
    quotient of eigenvalues of $\tilde{\omega }$ at $P$ is not in ${\mathbb Q}_{<0}$.    
  \end{enumerate}
\end{lemma}

\section{Rank $1$ valuations and transcendence of separatrices}
As valuations of ${\mathcal M}$ of rank one and rational rank one are exactly the divisorial ones,
those of contact with a non-convergent  curve and those having an infinite number of
Puiseux pairs, we have the following corollary, concerning the
algebraic independence of solutions of differential equations:
\begin{corollary}
  \label{cor:trascendence}
  Keep the notation of the previous sections.
  \begin{itemize}
  \item[a)] There is a valuation $\nu$ of ${\mathcal M}$, of rank $1$ and rational rank $1$ which is a
    \lh valuation for a vector field $\partial $ if and only if the foliation defined by $\partial $ is
    either dicritical or admits a separatrix $\hat{f}\in \hat{{\mathcal O}}$ which is 
    transcendental  over ${\mathcal O}$: that is, $(\hat{f})\cap {\mathcal O}=(0)$.
  \item[b)] If $\partial $ is non-dicritical, the existence of a rank $1$ and rational rank $1$ \lh
    valuation for $\partial $ is equivalent to $\partial $ having a transcendental  (that is,
    non-convergent) separatrix.
  \end{itemize}
\end{corollary}
More specifically, for vector fields with polynomial coefficients, by imitating the same
arguments as in Theorem \ref{the:formal-separatrices} and Corollary
\ref{cor:trascendence}, one can prove:
\begin{corollary}
  \label{cor:trascendence-polyn}
  Let $\partial =a(x,y)\partial/\partial x + b(x,y)\partial /\partial y$ be a polynomial vector
  field. Then
  \begin{itemize}
  \item[a)] There is a valuation $\nu$ of ${\mathbb C}(x,y)$, of rank $1$ and rational rank $1$ which is
    a \lh valuation for  $\partial $ if and only if the foliation defined by $\partial $ is
    either dicritical or admits a separatrix $\hat{f}\in {\mathbb C}[[x,y]]$ which is algebraically
    transcendental  over ${\mathbb C}(x,y)$.
  \item[b)] If $\partial $ is non-dicritical, the existence of a rank $1$ and rational rank $1$ \lh
    valuation for $\partial $ is equivalent to $\partial $ having an algebraically transcendental  separatrix.
  \item[c)] Moreover, if $\partial $ is regular, the (unique analytic) solution of $\partial $ passing
    through $(0,0)$ is transcendental  over ${\mathbb C}(x,y)$ if and only if there is a 
    \lh  valuation 
    for $\partial $ in ${\mathbb C}(x,y)$ with rank $1$ and rational rank $1$.
  \end{itemize}
\end{corollary}

\section{The ``missing'' separatrices}
In the previous sections, we have characterized \emph{some} of the separatrices of a germ
of analytic vector field in terms of their associated valuations: those having generic or
infinite character. The remaining cases include the solution of regular analytic vector
fields ---which have always rational character. 
Recall that in the proof of Theorem \ref{the:generic-infinite}, in
order  to show that
separatrices passing through singularities with rational quotient of eigenvalues did not give rise to \lh
valuations, we chose the functions:
\begin{equation*}
  a=u^{m}v^{n+1},\text{ and } b=u^{m}v^{n},
\end{equation*}
$(u=0)$ being the equation of the separatrix. In fact, any counterexample to \lhs property must be of
that kind: it has to include as a factor, the equation of the separatrix under study.  This led us to
give the following
\begin{definition}
  \label{def:weakly-lhopital}
  Keep the notation of the previous sections. We say that a valuation $\nu:{\mathcal
  M}^{\star}\rightarrow 
  \Gamma$ is a {\em weakly \lh valuation} for a vector field $\partial $ if 
  \begin{equation*}
    \nu \left( \frac{a}{b}-\frac{\partial a}{\partial b} \right)>0
  \end{equation*}
  for any $a,b\in{\mathcal M}$
  with $\nu(a)\geq \nu(b)>0$ {\em and such that if $\nu(c)>0$ then there is $n\in {\mathbb N}$ with
  $n\nu(c)>\nu(a)$}.
\end{definition}
This property can be expressed using the notion of {\em isolated subgroups}:
\begin{definition*}
  Let $ \Gamma$ be an ordered group. An {\em isolated subgroup of $\Gamma$} is a subgroup
  $(0)\neq\text{B}\subset \Gamma$ which is a segment: this means
  that, if $a,b\in \text{B}$ and $a\leq c\leq b$, then $c\in \text{B}$.
\end{definition*}
Isolated subgroups are ordered by inclusion, and, its number coincides with the rational rank
of the valuation, in our context. The {\em first} isolated subgroup is the smallest one. One can see that the first
isolated subgroup $\Gamma_{0}$ of a group $\Gamma$ is always $\Gamma-$ archimedian, which means that
given $a\in \Gamma_{0}$ and $\gamma \in \Gamma$ with $0<\gamma<a$, there is $n\in {\mathbb N}$ such that
$n\gamma >a$. In other words, the first isolated subgroup is the biggest subgroup ``not having infinitely
great elements''. Thus, Definition \ref{def:weakly-lhopital} becomes
\begin{proposition*}
  \label{def:weakly-lhopital2}
  A valuation $\nu$ is \emph{weakly \lh} for $\partial $ if and only if
  \begin{equation*}
    \nu \left( \frac{a}{b}-\frac{\partial a}{\partial b} \right)>0,    
  \end{equation*}
  for any $a,b$ such that $\nu(a),\nu(b)$ are in the first isolated subgroup of $\Gamma$ and 
  with $\nu(a)\geq\nu(b)>0$.
\end{proposition*}
The proof of the following theorem is straightforward
\begin{theorem}
  \label{the:weak-separatrices}
  Let $\nu$ be a valuation centered in ${\mathcal O}$ associated to a curve $f\in{\mathcal O}$ and let
  $\partial $  be a
  (germ of) vector field at $({\mathbb C}^{2},0)$. Then the following conditions are equivalent
  \begin{enumerate}
  \item $(f=0)$ is a separatrix of $\partial $.
  \item $\nu$ is a weakly \lh valuation for $\partial $.
  \end{enumerate}
\end{theorem}
From Theorem \ref{the:weak-separatrices} we infer the ``valuative Cauchy Theorem'':
\begin{theorem}
  \label{the:valuative-cauchy}
  Let $\partial $ be a regular (germ of) vector field at $({\mathbb C}^{2},0)$. There is one and only one
  irreducible (principal) ideal $(f)\in{\mathcal O}$ such that its associated valuation $\nu_{f}$ is a
  \lh valuation for $\partial $.
\end{theorem}
\providecommand{\bysame}{\leavevmode\hbox to3em{\hrulefill}\thinspace}

\end{document}